\documentclass[11pt]{article}
\font\smallit=cmti12
\usepackage{amssymb,amsmath,latexsym,amsthm}

\usepackage[linesnumbered, lined,boxed, commentsnumbered]{algorithm2e}

\newenvironment{frcseries}{\fontfamily{frc}\selectfont}{}

\newcommand\md[1]{\,\,(\mathrm{mod }\, #1)}
\def\s{S_{\mathrm{\mathfrak{z}}}}
\def\stwo{S_{\mathrm{\mathfrak{z,2}}}}

\def\B{\hfill $\Box$ \vskip 5pt}
\def\E{$\mathcal{E}$ }
\def\EE{$\mathcal{E}$}
\def\dd{\hfill $\diamond$ \vskip 5pt}
\def\d{\hfill $\dag$ \vskip 5pt}

\newtheorem{thm}{Theorem}

\newtheorem{prop}[thm]{Proposition}
\theoremstyle{definition}
\newtheorem{defn}[thm]{Definition}

\makeatletter

\renewcommand\section{\@startsection {section}{1}{\z@}%
{-30pt \@plus -1ex \@minus -.2ex}%
{2.3ex \@plus.2ex}%
{\normalfont\normalsize\bfseries}}

\renewcommand\subsection{\@startsection{subsection}{2}{\z@}%
{-3.25ex\@plus -1ex \@minus -.2ex}%
{1.5ex \@plus .2ex}%
{\normalfont\normalsize\bfseries}}
\renewcommand{\@seccntformat}[1]{\csname the#1\endcsname. } %\quad}

\makeatother

\begin{document}

\begin{center}
\uppercase{\bf  Zero-sum Generalized Schur Numbers}\\[20pt]
{\bf Aaron Robertson} \\
{\smallit Department of Mathematics,
Colgate University, Hamilton, New York}\\
{\tt arobertson@colgate.edu}
\end{center}

\vskip 10pt
\centerline{\bf Abstract}
\noindent
Let $r$ and $k$ be  positive integers
with $r \mid k$.  Denote by $\s(k;r)$ the minimum integer $n$ such
that every coloring $\chi:[1,n] \rightarrow \{0,1,\dots,r-1\}$ admits a solution to $\sum_{i=1}^{k-1} x_i = x_k$
with $\sum_{i=1}^{k} \chi(x_i) \equiv 0 \md{r}$. 
We give some formulas and lower bounds for various instances.

\baselineskip=14pt

\section{Introduction}

We start with the definition of the standard generalized Schur numbers.  For any
positive integers $k$ and $r$, there exists a minimal integer $S(k;r)$ such that
any $r$-coloring $\chi$ of $[1,S(k;r)]$ admits a monochromatic solution to
$\sum_{i=1}^{k-1}x_i = x_k$.  This follows directly from Ramsey's theorem
by defining the coloring of each edge $ij$ of the complete graph $K_n$
to be $\chi(|j-i|)$.  A monochromatic $K_k$ under this coloring with
vertices $i_1<i_2<\dots<i_k$ means that $i_{j+1}-i_j$ and $i_k-i_1$ are all
the same color under $\chi$.  Letting $x_j = i_{j+1}-i_j$ for $j=1,2,\dots,k-1$
and $x_k=i_k-i_1$ finishes the proof.  This is a generalization of the Schur
numbers, which are the special case $k=3$.  (Note that these definitions do
not agree with those found in \cite{BB}.)
An alternative method of showing that $S(k;r)$ exists would be to provide an
upper bound for it.  Beutelspacher and   Brestovansky \cite{BB} showed that
$S(k;2) = k^2-k-1$ thereby providing the independent existence of $S(k;r)$
for $r=2$.  

In this article we change the monochromatic property to a zero-sum property.  

\begin{defn} Let $a_1,a_2,\dots,a_n$ be a sequence of non-negative integers
and let $m \in \mathbb{Z}^+$.  We say that the sequence is {\it $m$-zero-sum} if
$\sum_{i=1}^n a_i \equiv 0 \md{m}$.
\end{defn}

The foundational zero-sum result  is the
Erd\H{o}s-Ginzberg-Ziv theorem \cite{EGZ}, which states that
any sequence of $2n-1$ integers contains an $n$-zero-sum subsequence of
$n$ integers.  Since around 1990, research activity concerning zero-sum
results has flourished,   through both the lens of additive number theory and
Ramsey theory.  An important extension of the Erd\H{o}s-Ginzberg-Ziv theorem
is
the weighted Erd\H{o}s-Ginzberg-Ziv theorem due to Grynkiewicz \cite{G}. It allows us to multiply the integers
in the Erd\H{o}s-Ginzberg-Ziv theorem by weights;
in particular,   if $w_1,w_2,\dots,w_n$ is an $n$-zero-sum sequence and
$a_1,a_2,\dots,a_{2n-1}$ is a sequence of $2n-1$ integers, then there exists
an $n$-term subsequence  $a_{i_1}, a_{i_2},\dots,a_{i_n}$ and a permutation
$\pi$ of $\{i_1,i_2,\dots,i_n\}$ such that
$\sum_{j=1}^n w_j a_{\pi(i_{j})} \equiv 0\md{n}$.
Further recent results can be found in  \cite{AM}, \cite{BCRY}, and \cite{GG}
among many others.

Most investigations of zero-sum sequences do not have a structure imposed on them.  This is in contrast to zero-sum results on edgewise colored graphs,
which have been around for many years (see, e.g., \cite{AC}, \cite{B}, \cite{BD}, and \cite{C} ).
Some notable exceptions are found in works of Bialostocki, such as \cite{BBS} and \cite{BSS} where the zero-sum sequence $x_1,x_2,\dots,x_n$ satisfies
$\sum_{i=1}^{n-1} x_i < x_n$ and in \cite{BBCY} where $x_{i+1}-x_{i} \leq x_i-x_{i-1}$ for $1 \leq i \leq n-1$.
These exceptions, however, do not have a rigid structure imposed on them due to the use of
inequality.
Very recently, a rigid structure similar to what we are investigating in this article was
investigated in \cite{Br}, while
in \cite{R}, the current author investigated 
zero-sum arithmetic progressions.  This article 
continues investigation of zero-sum sequences with an imposed rigid structure.

Throughout the paper we let $\mathcal{E}$ represent the equation $\sum_{i=1}^{k-1} x_i=x_k$.

\begin{defn} Let $k,r \in \mathbb{Z}^+$ with
$r \mid k$. We denote by $\s(k;r)$ the minimum integer such
that every coloring of $[1,\s(k;r)]$ with the colors $0,1,\dots,r-1$  admits an   $r$-zero-sum solution to \E.
We denote by $\stwo(k;r)$ the minimum integer such
that every coloring of $[1,\stwo(k;r)]$ with the colors $0$ and $1$ admits an $r$-zero-sum solution to \EE.
\end{defn}

The above definition assumes the existence of the respective minimum numbers.  Existence follows directly from the existence of the generalized Schur numbers $S(k;r)$.
  Note that
we need only prove the existence of $\s(k;r)$ since we easily
have $\stwo(k;r) \leq \s(k;r)$ as $\mathbb{Z}_2 \subseteq \mathbb{Z}_r$.  The existence of $\s(k;r)$ comes
from $\s(k;r) \leq S(k;r)$ as any  monochromatic solution to \E
 is $r$-zero-sum when $r \mid k$.  When $r \nmid k$, 
coloring every integer of $\mathbb{Z}^+$ with the color $1$ does not
admit a $k$-term $r$-zero-sum solution to \E and so we
write $\s(k;r)=\stwo(k;r)=\infty$ in this situation.

\section{Some Calculations}

The author wrote the fortran programs {\tt ZSGS.f} and {\tt ZSGS2.f},
available  at {\tt www.aaronrobertson.org}, to determine the
numbers $\s(k;r)$ and $\stwo(k;r)$, respectively, for small values of $k$ and $r$.
In addition to a standard backtrack algorithm for traversing colorings, we must
have a quick subroutine to determine solutions to $\sum_{i=1}^{k-1}x_i = x_k$
since checking $\approx n^k$ possible arrays $(x_1,x_2,\dots,x_k)$
on $[1,n]$ will quickly
become problematic.
To this end, in Algorithm 1, below, we give the pseudocode for our recursive subroutine.
In the code, we assume $x_1 \leq x_2 \leq \cdots \leq x_{k-1}$.

\IncMargin{1em}
\begin{algorithm}[h!] \small
\BlankLine
\SetKwFunction{solutions}{solutions}
\SetKwProg{myalg}{Subroutine}{}{}
\SetKwInOut{Input}{inputs}\SetKwInOut{Output}{output}
\Input{$t$, $k$\;}
\Output{Set $S$ of  solutions to $\sum_{i=1}^{k-1}x_i = t$\;}
Let $n=k-1$\;
Let $R$ be a partial solution, initialized to the array $(-,-,\dots,-,t)$ of length $k$\;
Let $S$ be the set of solutions, initialized empty\;
{\bf call} {\tt solutions(}$t,n,S,R${\tt )}\;
{\bf return}\;
\BlankLine
\BlankLine

\myalg{\solutions{$t,n,S,R$}}
{
\If{$t \leq 0$}{{\bf return}\;}

\eIf{$n = 1$}{
We have found a solution ($R$ is complete) so we add $R$ to $S$ and
{\bf return}\;
}
{
\For{$i$ from $\lceil \frac{t}{n}\rceil$ to $R_{n+1}$}{
$R_n=i$\;
{\bf call} {\tt solutions(}$t-i,n-1,S,R${\tt )}\;
}
}
}

\KwRet   $S$\;
\BlankLine
\caption{Solutions to $\sum_{i=1}^{k-1} x_i = t$}
\end{algorithm}

\newpage
Using this algorithm along with standard backtracking, we calculated the following values.

\begin{center}
\begin{tabular}{||l|c|c|c|c||} \hline\hline
\multicolumn{1}{||c|}{${}_{\mbox{\large $k$}} \,\,\,\,{}_\diagdown {\,\,\mbox{\large $r$}}$}      &       \multicolumn{1}{c|}{2}
   &  \multicolumn{1}{c|}{3}  &  \multicolumn{1}{c|}{4}  &
    \multicolumn{1}{c||}{5}  \\[5pt] \hline
    $2$                   &  1  &  $\infty$  &  $\infty$  &  $\infty$      \\[-1pt] 
    $3$                   &  $\infty$   &  $10$  &  $\infty$  &  $\infty$        \\[-1pt] 
    $4$                   &  5  &  $\infty$  &  $13$  &  $\infty$    \\[-1pt] 
    $5$                   &  $\infty$  &  $\infty$  &   $\infty$ &$38$     \\[-1pt] 
    $6$                   &  9  &  $15$  &  $\infty$  &  $\infty$      \\[-1pt] 
    $7$                   &  $\infty$  &  $\infty$  &  $\infty$  &  $\infty$    \\[-1pt] 
    $8$                   &  13  &  $\infty$  &  $27$  &  $\infty$     \\[-1pt] 
    $9$                   &  $\infty$  &  $24$  &  $\infty$  &  $\infty$       \\[-1pt] 
    $10$                  &  17  &  $\infty$  &  $\infty$  &  $\geq 45$      \\[-1pt] 
    $11$                  &  $\infty$  &  $\infty$  &  $\infty$  &  $\infty$      \\[-1pt] 
    $12$                  &  21  &  $33$  &  $43$  &  $\infty$     \\
\hline\hline
    \end{tabular}
\vspace*{-5pt}
\begin{center}
{\small {\bf Table 1}:  Values and a lower bound for $\s(k;r)$ for small $k$ and $r$.\\  The lower bound was the best
one achieved after 28 days of computing time.}
\end{center}
    \end{center}

\vskip 10pt

\begin{center}
\begin{tabular}{||l|c|c|c|c||} \hline\hline
\multicolumn{1}{||c|}{${}_{\mbox{\large $k$}} \,\,\,\,{}_\diagdown {\,\,\mbox{\large $r$}}$}      &       \multicolumn{1}{c|}{2}
   &  \multicolumn{1}{c|}{3}  &  \multicolumn{1}{c|}{4}  &
    \multicolumn{1}{c||}{5}  \\[5pt] \hline
    $2$                   &  1  &  $\infty$  &  $\infty$  &  $\infty$      \\[-1pt] 
    $3$                   &  $\infty$   &  $5$  &  $\infty$  &  $\infty$        \\[-1pt] 
    $4$                   &  5  &  $\infty$  &  $11$  &  $\infty$    \\[-1pt] 
    $5$                   &  $\infty$  &  $\infty$  &   $\infty$ &$19$     \\[-1pt] 
    $6$                   &  9  &  $13$  &  $\infty$  &  $\infty$      \\[-1pt] 
    $7$                   &  $\infty$  &  $\infty$  &  $\infty$  &  $\infty$    \\[-1pt] 
    $8$                   &  13  &  $\infty$  &  $25$  &  $\infty$     \\[-1pt] 
    $9$                   &  $\infty$  &  $22$  &  $\infty$  &  $\infty$       \\[-1pt] 
    $10$                  &  17  &  $\infty$  &  $\infty$  &  $41$      \\[-1pt] 
    $11$                  &  $\infty$  &  $\infty$  &  $\infty$  &  $\infty$      \\[-1pt] 
    $12$                  &  21  &  $31$  &  $41$  &  $\infty$     \\
\hline\hline
    \end{tabular}
\vspace*{-5pt}
\begin{center}
{{\bf Table 2}:  Values  for $\stwo(k;r)$ for small $k$ and $r$}
\end{center}
    \end{center}

\section{Some Formulas and Bounds}

\begin{prop} Let $k$ be an even positive integer.  Then $\s(k;2)=\stwo(k;2)=2k-3$.
\end{prop}

\noindent
{\it Proof.}
The fact that $\s(k;2)=\stwo(k;2)$ is by definition so we need only show
that $\s(k;2)=2k-3$. The formula obviously holds for $k=2$ so
we may assume that $k \geq 4$. To see that $\s(k;2) \geq 2k-3$ consider the $2$-coloring of $[1,2k-4]$
defined by coloring every integer in $[1,k-2]$ with color $0$ and every integer in
$[k-1,2k-4]$ with color $1$.  (In the sequel, we will describe this coloring by
$0^{k-2}1^{k-2}$.) If $x_1,x_2,\dots,x_{k-1}$ are all of color $0$, then
$x_k$ must be of color $1$ since $\sum_{i=1}^{k-1} x_i \geq k-1$.  Assuming that
$x_1 \leq x_2 \leq \cdots \leq x_{k-1}$, we can assume that $x_{k-1}$ has color $1$.
But then $\sum_{i=1}^{k-1} x_i \geq 2k-3$, a contradiction.

Moving on to the upper bound, consider an arbitrary coloring $\chi: [1,2k-3] \rightarrow \{0,1\}$.
Assume for a contradiction that $\chi$ does not admit a $2$-zero-sum solution to \EE.
We may assume $\chi(1)=0$ since $\chi$ admits a 2-zero-sum solution if and only if
$\widehat{\chi}$ defined by $\widehat{\chi}(i)=1-\chi(i)$ does.
Using $\chi(1)=0$ we deduce that
$\chi(k-1)=1$.  Considering the solution
$(x_1,x_2,\dots,x_k) = (1,1,\dots,1,k-1,2k-3)$ we conclude that $\chi(2k-3)=0$.
We deduce that $\chi(2)=1$ by considering the solution $(1,2,2,\dots,2,2k-3)$.
Consequently, we have $\chi(k)=0$ since $1+1+\cdots+1+2=k$.  Similarly,
by considering $1+1+\cdots+1+2+2=k+1$ we have $\chi(k+1)=1$.  Now, if $k=4$ then
$k+1=2k-3$ and we have a contradiction, so we may assume that $k \geq 6$.

By considering the solution $(1,1,\dots,1,2,k-2,2k-3)$ we deduce that $\chi(k-2)=0$.
This implies that $\chi(2k-4)=1$ since $1+1+\dots+1+(k-2)=2k-4$.
By considering the solution $(1,1,\dots,1,3,k)$ we conclude that $\chi(3)=1$.
Hence, we find that $\chi(k-4)=1$ since $1+1+\cdots+1+3+(k-4)=2k-4$.
We finish by noting that $(1,1,\dots,1,2,2,k-4,2k-4)$ is a
$2$-zero-sum solution to \EE, a contradiction.
\B

\begin{thm} Let $k \in \mathbb{Z}^+$ with $3 \mid k$.  Then $\s(k;3) \geq 3k-3$.
\label{th4}
\end{thm}

\noindent
{\it Proof.} We prove this by giving a coloring $\chi:[1,3k-4] \rightarrow \{0,1,2\}$
that avoids $3$-zero-sum solutions to \EE.  To this end, define $\chi$ by
$$
(012)^{\frac{k}{3}-1} (011)^{\frac{k}{3}} (021)^{\frac{k}{3}-1} 02.
$$
We will show that no solution to \E is 3-zero-sum under $\chi$.  We assume that $x_1 \leq x_2 \leq \cdots
\leq x_{k-1}$ and will use the notation
$$
T(m) = \sum_{i=1}^{m} x_i \qquad \mbox{and} \qquad C(m) = \sum_{i=1}^{m} \chi(x_i). 
$$

For an arbitrary solution to \E given by $(x_1,x_2,\dots,x_k)$, we let $A_j$ be the set of $x_i$ of color $j$, for $j=0,1,2$ restricted to $i \leq k-1$.
\vskip 5pt
\noindent
{\tt Case I.} $x_{k-1} \leq k-1$.  Since $x_{k-1} \leq k-1$, we see that, for  $1 \leq i \leq k-1$, 
in this case we have:

\vskip 5pt \hskip 114pt
(1) $\chi(x_i)= 0$ if and only if $x_i \equiv 1 \md{3}$;
\vskip 5pt\hskip 114pt
(2) $\chi(x_i)= 1$ if and only if $x_i \equiv 2 \md{3}$; and
\vskip 5pt\hskip 114pt
(3) $\chi(x_i)= 2$ if and only if $x_i \equiv 0 \md{3}$.

\vskip 5pt
\noindent
{\tt Subcase i.} $T(k-1) \equiv 0 \md{3}$. We must have $\chi(x_k)=1$ since $x_k \equiv 0 \md{3}$
and $x_k \geq k-1$.  We will show that   $C(k-1) \not\equiv 2 \md{3}$ thereby showing that there is no $3$-zero-sum solution
to \E in this subcase.  Assume, for a contradiction, that there exists a solution with
$C(k-1) \equiv 2 \md{3}$.  

We know that
$T(k-1) \equiv |A_0|+2|A_1|\equiv (k-1-|A_1|-|A_2|)+2|A_1| \equiv |A_1| - |A_2| + k-1 \equiv |A_1| - |A_2| + 2 \md{3}$.
Hence, $|A_2| \equiv |A_1|+2 \md{3}$.  
Using this, we have $C(k-1) \equiv  |A_1|+2|A_2| \equiv 3|A_1|+4
\equiv 1\md{3}$, contradicting our assumption that $C(k-1) \equiv 2 \md{3}$.
\d

\noindent
{\tt Subcase ii.} $T(k-1) \equiv 1 \md{3}$. We must have $\chi(x_k)=0$ since
$x_k \equiv 1 \md{3}$. We will show that we do not have $C(k-1) \equiv 0 \md{3}$. 
Following the argument in Subcase i, we have
$T(k-1) \equiv |A_1|-|A_2|+2 \md{3}$ so that $|A_2| \equiv |A_1|+1\md{3}$.  Then we
have $C(k-1) \equiv |A_1|+2|A_2| \equiv 3|A_1|+ 2$ and we conclude that $C(k-1) \equiv 2 \md{3}$
so that our solution is not 3-zero-sum.\d

\noindent
{\tt Subcase iii.} $T(k-1) \equiv 2 \md{3}$. We must have $\chi(x_k)\neq 0$ since
$x_k \equiv 2 \md{3}$. 
Following the argument in Subcase i, we have
$T(k-1) \equiv |A_1|-|A_2|+2 \md{3}$ so that $|A_2| \equiv |A_1|\md{3}$.  Then we
have $C(k-1) \equiv |A_1|+2|A_2| \equiv 3|A_1|$ and we conclude that $C(k-1) \equiv 0 \md{3}$.
Since $\chi(x_k) \neq 0 \md{3}$,
our solution is not 3-zero-sum.\d

This completes Case I. \dd

\vskip 5pt
\noindent
{\tt Case II.} $x_{k-1} \geq k$. In order to have $x_k \leq 3k-4$ we must have $x_{k-2} \leq k-1$.
Since $x_{k-2} \leq k-1$, we can use the arguments in Case I by considering $A_j$ restricted to $i \leq k-2$.
To this end, let $B_j = A_j \setminus \{x_{k-1}\}$ for $j=0,1,2$.

 \vskip 5pt
\noindent
{\tt Subcase i.} $T(k-2) \equiv 0 \md{3}$. Notice that we must have $x_{k-1} \equiv x_k \md{3}$
in this subcase. Hence,
we know that $\chi(x_{k-1})+\chi(x_k) \not\equiv 1 \md{3}$ (we cannot have $\chi(x_{k-1})=2$ 
since this gives $x_k > 3k-4$, which is out of bounds).
We will show that we must have $C(k-2) \equiv 2 \md{3}$ so that we cannot
have $C(k) \equiv 0 \md{3}$.  Using the arguments in Case I  we can conclude
that $|B_2| \equiv |B_1|+1 \md{3}$. From here we deduce that
$C(k-2) \equiv 2 \md{3}$, so that $C(k) \not\equiv 0 \md{3}$, and we are done with this subcase. \d
 
\noindent
{\tt Subcase ii.} $T(k-2) \equiv 1 \md{3}$. In this situation we must have
$1+x_{k-1} \equiv x_k \md{3}$.  Looking at $\chi$ we see that  
 $\chi(x_{k-1}) + \chi(x_k) \not\equiv 0\md{3}$.  We will show that
$C(k-2) \equiv 0 \md{3}$ so that $\chi$ does not contain a 3-zero-sum
solution to \E in this subcase.  Using the arguments in Case I  we conclude that 
$T(k-2) \equiv |B_1| - |B_2|+1 \md{3}$ so that $|B_1| \equiv |B_2| \md{3}$.
This gives us $C(k-2) \equiv 3|B_1| \equiv 0 \md{3}$, finishing this subcase. \d

\noindent
{\tt Subcase iii.} $T(k-2) \equiv 2 \md{3}$. In this situation we must have
$\chi(x_{k-1})+\chi(x_k) \not\equiv 2 \md{3}$.  We will show that
$C(k-2) \equiv 1 \md{3}$ so that $\chi$ does not contain a 3-zero-sum
solution to \E in this subcase.  Using the arguments in Case I  we conclude that 
$T(k-2) \equiv |B_1| - |B_2|+1 \md{3}$ so that $|B_2| \equiv |B_1|+2 \md{3}$.
This gives us $C(k-2) \equiv 3|B_1|+4 \equiv 1 \md{3}$, finishing this subcase. \d

This concludes the proof of Case II. \dd

 Having exhausted all possibilities,
the proof is complete.
\B

When we restrict the number of colors to just two, we can provide a formula for the associated number.
\begin{thm} Let $k \in \mathbb{Z}^+$ with $3 \mid k$.  Then $\stwo(k;3)=3k-5$.
\end{thm}

\noindent
{\it Proof.} To see that $\stwo(k;3) > 3k-6$ consider the coloring $0^{k-2}1^{2k-4}$.
In any solution to \E we must have at least one integer of color 1.  In order
to be 3-zero-sum we must then have at least 3 integers of color 1.  But then
$\sum_{i=1}^{k-1} x_i \geq 1(k-3) + 2(k-1) = 3k-5 > 3k-6$ so we cannot have
a solution with more than 2 integers of color 1.  

To show that $\stwo(k;3) \leq 3k-5$, assume, for a contradiction, that $\chi:[1,3k-5] \rightarrow
\{0,1\}$ does not admit a 3-zero-sum solution to \EE.  We may assume that $\chi(1)=0$
by considering $\widehat{\chi}(i) = 1 - \chi(i)$ and noticing that a solution is 3-zero-sum under
$\chi$ if and only if the solution is 3-zero-sum under $\widehat{\chi}$ (by the divisibility
property of $k$).
Considering the solution $(1,1,\dots,1,k-1)$ we must have $\chi(k-1)=1$.

\vskip 5pt
\noindent
{\tt Case I.} $\chi(2)=0$. Since
$\chi(k-1)=1$, from $(1,1,\dots,1,k-1,k-1,3k-5)$
we see that $\chi(3k-5)=0$.  In turn, since $1+3+3+\cdots+3=3k-5$ we see that $\chi(3)=1$.
Finally, consider $(2,2,3,3,\dots,3,3k-5)$.  The sum of the colors for this solution
is $(k-3)$, which is congruent to 0 modulo 3 since $3 \mid k$, a contradiction. \d

\vskip 5pt
\noindent
{\tt Case II.} $\chi(2)=1$.  Since $2+2+\cdots+2 = 2k-2$, we have $\chi(2k-2)=0$.
We also have $\chi(3k-5)=0$ by considering the solution $(2,2,\dots,2,k-1,3k-5)$.
In turn, since $1+3+3+\cdots+3=3k-5$ we have $\chi(3)=1$.  Now, for $0 \leq i \leq k-3$,
by considering the solution $({2,2,\dots,2},
3,3,2k)$ we have $\chi(2k)=0$.  
Next, consider $(1,1,\dots,1,k,2k-2)$ to see that $\chi(k)=1$.
But now
$(1,1,\dots,1,2,2,k,2k)$ is a 3-zero-sum solution to \EE, a contradiction. \d

As the two cases cover all situations, the proof is complete.\B

\begin{thm}
Let $k \in \mathbb{Z}^+$ with $4 \mid k$.  Then $\s(k;4) \geq 4k-5$.  
\end{thm}

\noindent
{\it Proof.} We prove this by giving a coloring $\chi:[1,4k-6] \rightarrow \{0,1,2,3\}$
that avoids $4$-zero-sum solutions to \EE.  To this end, define $\chi$ by
$$
(0123)^{\frac{k}{4}-1} (0120) (0220)^{\frac{k}{2}-1} (3210)^{\frac{k}{4}-1} 32.
$$
We will show that no solution to \E is 4-zero-sum under $\chi$.  We assume that $x_1 \leq x_2 \leq \cdots
\leq x_{k-1}$ and will again use the notation
$$
T(m) = \sum_{i=1}^{m} x_i \qquad \mbox{and} \qquad C(m) = \sum_{i=1}^{m} \chi(x_i). 
$$

For an arbitrary solution to \E given by $(x_1,x_2,\dots,x_k)$, we let $A_j$ be the set of $x_i$ of color $j$, for $j=0,1,2,3$ restricted to $i \leq k-1$.

\vskip 5pt
\noindent
{\tt Case I.} $x_{k-1} \leq k-1$.  Since $x_{k-1} \leq k-1$, we see that, for  $1 \leq i \leq k-1$, 
in this case we have:

\vskip 5pt \hskip 114pt
(1) $\chi(x_i)= 0$ if and only if $x_i \equiv 1 \md{4}$;
\vskip 5pt\hskip 114pt
(2) $\chi(x_i)= 1$ if and only if $x_i \equiv 2 \md{4}$; 
\vskip 5pt\hskip 114pt
(3) $\chi(x_i)= 2$ if and only if $x_i \equiv 3 \md{4}$; and
\vskip 5pt\hskip 114pt
(3) $\chi(x_i)= 3$ if and only if $x_i \equiv 0 \md{4}$.

\vskip 5pt
\noindent
{\tt Subcase i.} $T(k-1) \equiv 0 \md{4}$. We must have $\chi(x_k)=0$ since $x_k \equiv 0 \md{4}$.  We will show that   $C(k-1) \not\equiv 0 \md{4}$ thereby showing that there is no $4$-zero-sum solution
to \E in this subcase.  

We have
$T(k-1) \equiv |A_0|+2|A_1|+3|A_2|\equiv (k-1-|A_1|-|A_2|-|A_3|)+2|A_1|+3|A_2| 
\equiv |A_1| + 2|A_2| -|A_3|+3 \md{4}$.
Since we have $T(k-1) \equiv 0 \md{4}$ in this subcase, we conclude that
$|A_3| \equiv |A_1|+2|A_2|+3 \md{4}$.  
Using this, we have $C(k-1) \equiv  |A_1|+2|A_2|+3|A_3| \equiv 4|A_1|+8|A_2|+9
\equiv 1\md{4}$.  Hence, $C(k) \equiv 1 \md{4}$ so that there is no $4$-zero-sum
solution in this subcase.
\d

\noindent
{\tt Subcase ii.} $T(k-1) \equiv 1 \md{4}$. We must have $\chi(x_k)=0$ or $3$ since
$x_k \equiv 1 \md{4}$. We will show that $C(k-1)  \equiv 2\md{4}$ so that
we know $C(k)  \not\equiv 0 \md{4}$, and hence we do not have a 4-zero-sum solution. 
Following the argument in Subcase i, we have
$|A_3| \equiv |A_1|+2|A_2|+2\md{4}$.  Then we
have $C(k-1) \equiv |A_1|+2|A_2| + 3|A_3| \equiv   4|A_1|+8|A_2|+6 \equiv 2 \md{4}$ and we conclude that $C(k-1) \equiv 2 \md{3}$
so that $C(k) \equiv 1$ or $2 \md{4}$ and our solution is not 4-zero-sum.\d

\noindent
{\tt Subcase iii.} $T(k-1) \equiv 2 \md{4}$. We must have $\chi(x_k)= 2$ since
$x_k \equiv 2 \md{4}$ and we cannot have $\chi(x_k)=1$ since this means
$x_k= k-2$ , which is not possible.  
Following the argument in Subcase i, we have
$|A_3| \equiv |A_1|+2|A_2|+1\md{4}$. Hence, $C(k-1) \equiv |A_1|+2|A_2| + 3|A_3| \equiv   4|A_1|+8|A_2|+3 \equiv 3 \md{4}$,
so that
$C(k) \equiv 1 \md{4}$ and our solution is not 4-zero-sum.\d

\noindent
{\tt Subcase iv.} $T(k-1) \equiv 3 \md{4}$. We must have $\chi(x_k)= 1$ or $2$ since
$x_k \equiv 3 \md{4}$.
Following the argument in Subcase i, we have
$|A_3| \equiv |A_1|+2|A_2|\md{4}$. Hence, $C(k-1) \equiv |A_1|+2|A_2| + 3|A_3| \equiv   4|A_1|+8|A_2|  \equiv 0 \md{4}$. This gives us that
$C(k) \equiv 1$ or $2 \md{4}$ so that our solution is not 4-zero-sum.\d

This completes Case I.\dd

\vskip 5pt
\noindent
{\tt Case II.} $x_{k-2} \leq k-1$ and $x_{k-1} \geq k$.  Since $x_{k-2} \leq k-1$, we can use the arguments in Case I by considering $A_j$ restricted to $i \leq k-2$. Thus, we
let $B_j = A_j \setminus \{x_{k-1}\}$ for $j=0,1,2,3$.

\vskip 5pt
\noindent
{\tt Subcase i.} $T(k-2) \equiv 0 \md{4}$. We must have
$x_{k-1} \equiv x_k \md{4}$.  Looking at $\chi$, we see that
 $\chi(x_{k-1})+\chi(x_k) \in \{0,3\}$.  We will show that   $C(k-2) \not\equiv 0,1 \md{4}$ thereby showing that there is no $4$-zero-sum solution
to \E in this subcase.  

We have
$T(k-2) \equiv |B_0|+2|B_1|+3|B_2|\equiv (k-2-|B_1|-|B_2|-|B_3|)+2|B_1|+3|B_2| 
\equiv |B_1| + 2|B_2| -|B_3|+2 \md{4}$.
Since we have $T(k-2) \equiv 0 \md{4}$ in this subcase, we conclude that
$|B_3| \equiv |B_1|+2|B_2|+2 \md{4}$.  
Using this, we have $C(k-2) \equiv  |B_1|+2|B_2|+3|B_3| \equiv 4|B_1|+8|B_2|+6
\equiv 2\md{4}$, hence $C(k-2) \not\equiv 0,1 \md{4}$ and this subcase is complete.
\d

\noindent
{\tt Subcase ii.} $T(k-2) \equiv 1 \md{4}$. We must have $x_k\equiv x_{k-1}+1   \md{4}$.
From this we conclude that $\chi(x_{k-1})+\chi(x_k) \in \{0,2,3\}$.
We will show that $C(k-2)  \equiv 3\md{4}$ so that
we know $C(k)  \not\equiv 0 \md{4}$, and hence we do not have a 4-zero-sum solution. 
Following the argument in Subcase i, we have
$|B_3| \equiv |B_1|+2|B_2|+1\md{4}$.  Then we
have $C(k-2) \equiv |B_1|+2|B_2| + 3|B_3| \equiv   4|B_1|+8|B_2|+3 \equiv 3 \md{4}$ and we conclude that 
$C(k-2) \equiv 3 \md{3}$
so that $C(k) \equiv 1$ or $2 \md{4}$ and our solution is not 4-zero-sum.\d

\noindent
{\tt Subcase iii.} $T(k-2) \equiv 2 \md{4}$. 
Looking at our coloring, we see that we can only
have $\chi(x_{k-2})+\chi(x_{k-1})=0$ if $x_{k-1} \geq 3k-3$.
But then $T(k-1) > 4k-6$, which is out of bound.  Thus,
we must have $x_k\equiv x_{k-1}+2   \md{4}$.
From this we conclude that $\chi(x_{k-2})+\chi(x_{k-1}) \in \{1,2\}$.
Following the argument in Subcase i, we have
$|B_3| \equiv |B_1|+2|B_2|\md{4}$.  Then we
have $C(k-2) \equiv |B_1|+2|B_2| + 3|B_3| \equiv   4|B_1|+8|B_2| \equiv 0\md{4}$ and we conclude that 
$C(k) \equiv 1$ or $2 \md{4}$ so that our solution is not 4-zero-sum.\d

\noindent
{\tt Subcase iv.} $T(k-1) \equiv 3 \md{4}$. We must have $x_k\equiv x_{k-1}+3   \md{4}$.
As in Subcase iii directly above, we cannot have
$\chi(x_{k-2})+\chi(x_{k-1})=3$ since that would imply that $T(k-1)>4k-6$.
Hence, we see that $\chi(x_{k-2})+\chi(x_{k-1}) \in \{0,1,2\}$
Following the argument in Subcase i, we have
$|B_3| \equiv |B_1|+2|B_2|+3\md{4}$.  Then we
have $C(k-2) \equiv |B_1|+2|B_2| + 3|B_3| \equiv   4|B_1|+8|B_2|+9 \equiv 1\md{4}$
But then $C(k) \equiv 1,2,$ or $3 \md{4}$ so that $C(k) \equiv 0 \md{4}$ is not possible.

\vskip 5pt
This completes Case II.\dd

\vskip 5pt
\noindent
{\tt Case III.} $x_{k-2} \geq k$. We must have $x_{k-3} \leq k-1$ for otherwise
$T(k-1) \geq 4k-4$. Also, we have $x_k \geq 3k-3$ so there
is a one-to-one correspondence between $x_k$ and $\chi(x_k)$. Since $x_{k-3} \leq k-1$, we can use the arguments in Case I by considering $A_j$ restricted to $i \leq k-3$.
To this end, let $D_j = A_j \setminus \{x_{k-2}, x_{k-1}\}$ for $j=0,1,2,3$.
Note that, under $\chi$, the only possible colors of $x_{k-2}$ and $x_{k-1}$
are 0 and 2.

\vskip 5pt
\noindent
{\tt Subcase i.} $T(k-3) \equiv 0 \md{4}$.
We have
$T(k-3) \equiv |D_0|+2|D_1|+3|D_2|\equiv (k-3-|D_1|-|D_2|-|D_3|)+2|D_1|+3|D_2| 
\equiv |D_1| + 2|D_2| -|D_3|+1 \md{4}$.
Since we have $T(k-3) \equiv 0 \md{4}$ in this subcase, we conclude that
$|D_3| \equiv |D_1|+2|D_2|+1 \md{4}$.  
Using this, we have $C(k-3) \equiv  |D_1|+2|D_2|+3|D_3| \equiv 4|D_1|+8|D_2|+3
\equiv 3\md{4}$.

\vskip 5pt
\noindent
{\tt Subsubcase a.} $\chi(x_{k-2})=\chi(x_{k-1})=0$.  We have $x_{k-2}, x_{k-1} \equiv 0$ or $1\md{4}$
so that $T(k-1) \equiv 0,1,$ or $2\md{4}$.  We also know that $C(k-1) \equiv 3 \md{4}$ so we
need $\chi(x_k)=1$ in order to have a 4-zero-sum solution.  But the only possible integers of color 1
are congruent to $3$ modulo 4, which is not possible since $T(k-1) \not \equiv 3 \md{4}$.

\vskip 5pt
\noindent
{\tt Subsubcase b.} $\chi(x_{k-2})=0$ and $\chi(x_{k-1})=2$ or $\chi(x_{k-2})=2$ and $\chi(x_{k-1})=0$.
In order to have $C(k) \equiv 0 \md{4}$ we require $\chi(x_k)=3$, which yields $x_k \equiv 1 \md{4}$.
This means we need $x_{k-2}+x_{k-1} \equiv 1 \md{4}$ since we have $T(k-3) \equiv 0\md{4}$ in this subcase.
But we know that one of $x_{k-2}$ and $x_{k-1}$ has color 0, and so is congruent to $0$ or $1$
modulo 4, and the other has color 2, and so is congruent to $2$ or $3$  modulo 4.
Hence, we cannot have $x_{k-2}+x_{k-1} \equiv 1 \md{4}$ and, consequently, no 4-zero-sum solution to \E
exists.

\vskip 5pt
\noindent
{\tt Subsubcase c.} $\chi(x_{k-2})=\chi(x_{k-1})=2$. In order to have $C(k) \equiv 0 \md{4}$ we require 
$\chi(x_k)=1$, which yields $x_k \equiv 3 \md{4}$.  We know that $x_{k-2}$ and $x_{k-1}$ are
congruent to either $2$ or $3$ modulo 4 so that
$x_{k-2}+x_{k-1} \equiv 0,1,$ or $2\md{4}$.  But then $T(k-1) \not\equiv 3 \equiv x_k \md{4}$.  Hence, 
$C(k) \not \equiv 0 \md{4}$.

\vskip 5pt
This completes Subcase i. \d

\noindent
{\tt Subcase ii.} $T(k-3) \equiv 1 \md{4}$.
Following the argument in Case III.i we conclude that
$|D_3| \equiv |D_1|+2|D_2| \md{4}$ so that
$C(k-3) \equiv  |D_1|+2|D_2|+3|D_3| \equiv 4|A_1|+8|A_2|
\equiv 0\md{4}$.

\vskip 5pt
\noindent
{\tt Subsubcase a.} $\chi(x_{k-2})=\chi(x_{k-1})=0$.  We have $x_{k-2}, x_{k-1} \equiv 0$ or $1\md{4}$
so that $T(k-1) \equiv 1,2,$ or $3\md{4}$.  We also know that $C(k-1) \equiv 0 \md{4}$ so we
need $\chi(x_k)=0$ in order to have a 4-zero-sum solution.  But the only possible integers of color 0
are congruent to $0$ modulo 4, which is not possible since $T(k-1) \not \equiv 0 \md{4}$.

\vskip 5pt
\noindent
{\tt Subsubcase b.} $\chi(x_{k-2})=0$ and $\chi(x_{k-1})=2$ or $\chi(x_{k-2})=2$ and $\chi(x_{k-1})=0$.
In order to have $C(k) \equiv 0 \md{4}$ we require $\chi(x_k)=2$, which yields $x_k \equiv 2 \md{4}$.
This means we need $x_{k-2}+x_{k-1} \equiv 1 \md{4}$ since we have $T(k-3) \equiv 1\md{4}$ in this subcase.
But we know that one of $x_{k-2}$ and $x_{k-1}$ has color 0, and so is congruent to $0$ or $1$
modulo 4, and the other has color 2, and so is congruent to $2$ or $3$  modulo 4.
Hence, we cannot have $x_{k-2}+x_{k-1} \equiv 1 \md{4}$ and, consequently, no 4-zero-sum solution to \E
exists.

\vskip 5pt
\noindent
{\tt Subsubcase c.} $\chi(x_{k-2})=\chi(x_{k-1})=2$. In order to have $C(k) \equiv 0 \md{4}$ we require 
$\chi(x_k)=0$, which yields $x_k \equiv 0 \md{4}$.  We know that $x_{k-2}$ and $x_{k-1}$ are
congruent to either $2$ or $3$ modulo 4 so that
$x_{k-2}+x_{k-1} \equiv 0,1,$ or $2\md{4}$.  But then $T(k-1)  \equiv 1,2$ or $3\md{4}$
so that $T(k-1) \not \equiv x_k \md{4}$.  Hence, 
$C(k) \not \equiv 0 \md{4}$.

\vskip 5pt
This completes Subcase ii. \d

\noindent
{\tt Subcase iii.} $T(k-3) \equiv 2 \md{4}$.
Following the argument in Case III.i we   conclude that
$|D_3| \equiv |D_1|+2|D_2|+3 \md{4}$.  Hence,
$C(k-3) \equiv  |D_1|+2|D_2|+3|D_3| \equiv 4|A_1|+8|A_2|+9
\equiv 1\md{4}$.

\vskip 5pt
\noindent
{\tt Subsubcase a.} $\chi(x_{k-2})=\chi(x_{k-1})=0$.  We have $x_{k-2}, x_{k-1} \equiv 0$ or $1\md{4}$
so that $T(k-1) \equiv 0,2,$ or $3\md{4}$.  We also know that $C(k-1) \equiv 1 \md{4}$ so we
need $\chi(x_k)=3$ in order to be 4-zero-sum.  But the only possible integers of color 3
are congruent to $1$ modulo 4, which is not possible since $T(k-1) \not \equiv 1 \md{4}$.

\vskip 5pt
\noindent
{\tt Subsubcase b.} $\chi(x_{k-2})=0$ and $\chi(x_{k-1})=2$ or $\chi(x_{k-2})=2$ and $\chi(x_{k-1})=0$.
In order to have $C(k) \equiv 0 \md{4}$ we require $\chi(x_k)=1$, so that $x_k \equiv 3 \md{4}$.
This means we need $x_{k-2}+x_{k-1} \equiv 1 \md{4}$ since we have $T(k-3) \equiv 2\md{4}$ in this subcase.
As in Case III.ii.b, we cannot have $x_{k-2}+x_{k-1} \equiv 1 \md{4}$ and, consequently, no 4-zero-sum solution to \E
exists.

\vskip 5pt
\noindent
{\tt Subsubcase c.} $\chi(x_{k-2})=\chi(x_{k-1})=2$. In order to have $C(k) \equiv 0 \md{4}$ we require 
$\chi(x_k)=3$, so that $x_k \equiv 1 \md{4}$.  We know that $x_{k-2}$ and $x_{k-1}$ are
congruent to either $2$ or $3$ modulo 4 so that
$x_{k-2}+x_{k-1} \equiv 0,1,$ or $2\md{4}$.  But then $T(k-1)  \equiv 0,2$ or $3\md{4}$
so that $T(k-1) \not \equiv x_k \md{4}$.  Hence, 
$C(k) \not \equiv 0 \md{4}$.

\vskip 5pt
This completes Subcase iii. \d
 
\noindent
{\tt Subcase iv.} $T(k-3) \equiv 3 \md{4}$.
Following the argument in Case III.i we   conclude that
$|D_3| \equiv |D_1|+2|D_2|+2 \md{4}$ so that
$C(k-3) \equiv  |D_1|+2|D_2|+3|D_3| \equiv 4|A_1|+8|A_2|+6
\equiv 2\md{4}$.

\vskip 5pt
\noindent
{\tt Subsubcase a.} $\chi(x_{k-2})=\chi(x_{k-1})=0$.  We have $x_{k-2}, x_{k-1} \equiv 0$ or $1\md{4}$
so that $T(k-1) \equiv 0,1,$ or $3\md{4}$.  We also know that $C(k-1) \equiv 2 \md{4}$ so we
need $\chi(x_k)=2$ in order to be 4-zero-sum.  But the only possible integers of color 2
are congruent to $2$ modulo 4, which is not possible since $T(k-1) \not \equiv 2 \md{4}$.

\vskip 5pt
\noindent
{\tt Subsubcase b.} $\chi(x_{k-2})=0$ and $\chi(x_{k-1})=2$, or $\chi(x_{k-2})=2$ and $\chi(x_{k-1})=0$.
In order to have $C(k) \equiv 0 \md{4}$ we require $\chi(x_k)=0$, which yields $x_k \equiv 0 \md{4}$.
This means we need $x_{k-2}+x_{k-1} \equiv 1 \md{4}$ since we have $T(k-3) \equiv 3\md{4}$ in this case.
As in Case III.ii.b, we cannot have $x_{k-2}+x_{k-1} \equiv 1 \md{4}$ and, consequently, no 4-zero-sum solution to \E
exists.

\vskip 5pt
\noindent
{\tt Subsubcase c.} $\chi(x_{k-2})=\chi(x_{k-1})=2$. In order to have $C(k) \equiv 0 \md{4}$ we require 
$\chi(x_k)=2$, which yields $x_k \equiv 2 \md{4}$.  As in Case III.ii.c, we have
$x_{k-2}+x_{k-1} \equiv 0,1,$ or $2\md{4}$.  But then $T(k-1)  \equiv 0,1$ or $3\md{4}$
so that $T(k-1) \not \equiv x_k \md{4}$.  Hence, 
$C(k) \not \equiv 0 \md{4}$. 

\vskip 5pt
This is the end of the proof of Subcase iv.\d
 
Having covered all possibilities with $T(k-1)\equiv 3\md{4}$, we are done with Case III.\dd

Having exhausted all cases, the theorem's proof is complete. \B

The last instance we investigate are those numbers along the diagonal.

\begin{prop} Let $k$ be an odd positive integer.  Then $\s(k;k) \geq 2(k^2-k-1)$.
\end{prop}

\noindent
{\it Proof.} We will show that the $k$-coloring $(01)^{k-2} (0 \,(k-1))^{(k-1)(k-2)} (01)^{k-2} 0$
avoids $k$-zero-sum solutions to \EE.  It is easy to check that we cannot have
a solution to \E with all integers of color 1 or
all integers of color $k-1$.  Hence, in order to have a $k$-zero-sum solution,
the number of integers colored 1 must equal the number of integers colored $k-1$.
Now, $k$ being odd implies that we have an odd number of integers of color $0$.
Next, we note that the only integers of color $0$ are odd, while the only
integers of color $1$ or $k-1$ are even.  By comparing the parities
of $\sum_{i=1}^{k-1} x_i$ and $x_k$, this cannot occur.  Hence, we cannot
have a $k$-zero-sum solution to \E under this coloring. \B

If we restrict to just two colors, then the associated number is the same as
$S(k;2)$ since any $k$-zero-sum solution to \E must necessarily be monochromatic.
In other word,
for $k \in \mathbb{Z}^+$ we have $\stwo(k;k)=k^2-k-1$.

\section{Conclusion and Open Questions}  

The area of inquiry of zero-sum sequences with rigid structure is ripe for future research.  What
can be said about zero-sum Rado numbers in addition to what is found in \cite{Br}?  What can we
say about zero-sum sequences $x_1<x_2<\dots<x_k$ with $x_{i+1}-x_i$ from a prescribed set
(say the (shifted) primes, powers of 2, etc)?

For specific questions related to this article, we ask the following.

\begin{itemize} \renewcommand{\itemsep}{2pt}

\item[Q1.] Is it true that $\s(k;3)=3k-3$ for $k \geq 6$? 

\item[Q2.] Prove or disprove: $\s(k;4) = 4k-5$ for $k \geq 8$.

\item[Q3.] What is the exact value of $\stwo(k;4)$?

\item[Q4.] Is it true that $\s(k;k)$ is of order $k^2$?

\end{itemize}

 \end{document}